\documentclass[11pt]{amsart}
\usepackage{amssymb,latexsym}

\newdimen\AAdi%
\newbox\AAbo%
\def\AAk#1#2{\s_etbox\AAbo=\hbox{#2}\AAdi=\wd\AAbo\kern#1\AAdi{}}%
\def\AAr#1#2#3{\s_etbox\AAbo=\hbox{#2}\AAdi=\ht\AAbo\raise#1\AAdi\hbox{#3}}%
\font\tenmsb=msbm10 at 12pt
\font\sevenmsb=msbm7 at 8pt
\font\fivemsb=msbm5 at 6pt
\newfam\msbfam
\textfont\msbfam=\tenmsb
\scriptfont\msbfam=\sevenmsb
\scriptscriptfont\msbfam=\fivemsb
\def\Bbb#1{{\tenmsb\fam\msbfam#1}}

\newcommand{\beq}{\begin{equation}}
\newcommand{\eeq}{\end{equation}}
\newcommand{\beqr}{\begin{eqnarray}}
\newcommand{\eeqr}{\end{eqnarray}}
\newcommand{\ba}{\begin{array}}
\newcommand{\ea}{\end{array}}

\begin{document}

\newtheorem{thm}{Theorem}
\newtheorem{lem}{Lemma}
\newtheorem{cor}{Corollary}
\newtheorem{rem}{Remark}
\newtheorem{pro}{Proposition}
\newtheorem{defi}{Definition}
\newtheorem{conj}[thm]{Conjecture}
\newcommand{\noi}{\noindent}
\newcommand{\dis}{\displaystyle}
\newcommand{\mint}{-\!\!\!\!\!\!\int}

\def \bx{\hspace{2.5mm}\rule{2.5mm}{2.5mm}} \def \vs{\vspace*{0.2cm}}
\def\hs{\hspace*{0.6cm}}
\def \ds{\displaystyle}
\def \p{\partial}
\def \O{\Omega}
\def \o{\omega}
\def \b{\beta}
\def \m{\mu}
\def \l{\lambda}
\def\L{\Lambda}
\def \ul{u_\lambda}
\def \D{\Delta}
\def \d{\delta}
\def \k{\kappa}
\def \s{\sigma}
\def \e{\varepsilon}
\def \a{\alpha}
\def \tf{\tilde{f}}
\def\cqfd{%
\mbox{ }%
\nolinebreak%
\hfill%
\rule{2mm} {2mm}%
\medbreak%
\par%
}
\def \pr {\noindent {\it Proof.} }
\def \rmk {\noindent {\it Remark} }
\def \esp {\hspace{4mm}}
\def \dsp {\hspace{2mm}}
\def \ssp {\hspace{1mm}}

\def \u{u_+^{p^*}}
\def \ui{(u_+)^{p^*+1}}
\def \ul{(u^k)_+^{p^*}}
\def \energy{\int_{\R^n}\u }
\def \sk{\s_k}
\def \mo{\mu_k}
\def\cal{\mathcal}
\def \I{{\cal I}}
\def \J{{\cal J}}
\def \K{{\cal K}}
\def \OM{\overline{M}}

\def\fk{{{\cal F}}_k}
\def\M1{{{\cal M}}_1}
\def\Fk{{\cal F}_k}
\def\Fl{{\cal F}_l}
\def\FF{\cal F}
\def\Gk{{\Gamma_k^+}}
\def\n{\nabla}
\def\uuu{{\n ^2 u+du\otimes du-\frac {|\n u|^2} 2 g_0+S_{g_0}}}
\def\uuug{{\n ^2 u+du\otimes du-\frac {|\n u|^2} 2 g+S_{g}}}
\def\sku{\sk\left(\uuu\right)}
\def\qed{\cqfd}
\def\vvv{{\frac{\n ^2 v} v -\frac {|\n v|^2} {2v^2} g_0+S_{g_0}}}
\def\vvs{{\frac{\n ^2 \tilde v} {\tilde v}
 -\frac {|\n \tilde v|^2} {2\tilde v^2} g_{S^n}+S_{g_{S^n}}}}
\def\skv{\sk\left(\vvv\right)}
\def\tr{\hbox{tr}}
\def\pO{\partial \Omega}
\def\dist{\hbox{dist}}
\def\RR{\Bbb R}\def\R{\Bbb R}
\def\C{\Bbb C}
\def\B{\Bbb B}
\def\N{\Bbb N}
\def\Q{\Bbb Q}
\def\Z{\Bbb Z}
\def\PP{\Bbb P}
\def\EE{\Bbb E}
\def\F{\Bbb F}
\def\G{\Bbb G}
\def\H{\Bbb H}
\def\SS{\Bbb S}\def\S{\Bbb S}

\def\lcf{{locally conformally flat} }

\def\circledwedge{\setbox0=\hbox{$\bigcirc$}\relax \mathbin {\hbox
to0pt{\raise.5pt\hbox to\wd0{\hfil $\wedge$\hfil}\hss}\box0 }}

\def\sss{\frac{\s_2}{\s_1}}

\date{}
\title[A new conformal invariant] {A new conformal invariant on 3-dimensional manifolds
 }

\author{Yuxin Ge}
\address{Laboratoire d'Analyse et de Math\'ematiques Appliqu\'ees,
CNRS UMR 8050,
D\'epartement de Math\'ematiques,
Universit\'e Paris Est-Cr\'eteil Val de Marne \\61 avenue du G\'en\'eral de Gaulle,
94010 Cr\'eteil Cedex, France}
\email{ge@univ-paris12.fr}
\author{Guofang Wang}\thanks{YG is supported by ANR  project
ANR-08-BLAN-0335-01. GW is partly supported by SFB/TR71 ``Geometric partial differential equations''  of DFG}
\address{ Albert-Ludwigs-Universit\"at Freiburg,
Mathematisches Institut,
Eckerstr. 1,
D-79104 Freiburg, Germany}
\email{guofang.wang@math.uni-freiburg.de}
\begin{abstract}By improving the analysis developed in the study of $\s_k$-Yamabe problem, we  prove in this paper that the  De Lellis-Topping
inequality
is true on 3-dimensional
Riemannian manifolds of nonnegative scalar curvature. More precisely,
  if $(M^3, g)$ is a 3-dimensional closed Riemannian manifold with non-negative  scalar curvature, then
\[
\int_M |Ric-\frac{\overline R } 3 g|^2 dv (g)\le 9\int_M |Ric-\frac{ R } 3 g|^2dv(g),
\]
where $\overline R=vol (g)^{-1} \int_M R dv(g)$ is the average of the scalar curvature $R$ of $g$. Equality holds if and only if $(M^3,g)$ is a space form. We in fact study the following new conformal invariant
\[
\ds \widetilde Y([g_0]):=\sup_{ g\in {\cal C}_1([g_0])}\frac {\ds vol(g)\int_M \s_2( g) dv( g)}
{\ds (\int_M \s_1( g) dv( g))^2},
\]
where ${\cal C}_1([g_0]):=\{g=e^{-2u}g_0\,|\, R>0\}$ and prove that
$\widetilde Y([g_0])\le 1/3$, which implies the above inequality.

\end{abstract}
 \subjclass{Primary 53C21; Secondary 53C20, 58E11}

\maketitle
\section{Introduction}
Very recently, De Lellis and Topping proved an interesting result about a generalization of Schur Lemma

\

\noindent{\bf Theorem A.} [Almost Schur Lemma \cite{DT}] \label{thmDT} {\it For $n\ge 3$, if $(M^n, g)$ is an  $n$-dimensional 
closed Riemannian manifold with non-negative Ricci tensor, then
\begin{equation}\label{eq1}
\int_M |Ric-\frac{\overline R } n g|^2 dv (g)\le \frac {n^2}{(n-2)^2}\int_M |Ric-\frac{ R } n g|^2dv(g),
\end{equation}
where $\overline R=vol (g)^{-1} \int_M R dv(g)$ is the average of the scalar curvature $R$ of $g$.}

\

The result can be seen as a quantitative version
or a stability result of the Schur Lemma.
It was proved in \cite{DT} that  the constant in inequality (\ref{eq1}) is optimal and
the non-negativity of the Ricci tensor can not be removed in general:
 When $n\ge 5$  there are examples of metrics on $\S^n$
which make the radio  of the left hand side of (\ref{eq1}) to the
right hand side of (\ref{eq1}) arbitrarily large. When $n=3$, they
found manifolds which makes the ratio arbitrarily large. An
interesting question remains open: {\it Inequalities of this  form
may hold for n = 3 and n = 4 with constants depending on the
topology of M.}

With an observation that the De Lellis-Topping inequality is equivalent to an inequality in terms of
$\s_k$-scalar curvature 
\begin{equation}\label{eq3}
\left(\int_M \s_1(g) dv(g)\right)^2\ge \frac {2n}{n-1} vol(g) \int_M\s_2(g) dv(g),
\end{equation}
we proved in \cite{GeWang_Proc} that \eqref{eq1} holds for 4-dimensional manifolds
of nonnegative scalar curvature, by using an argument of Gursky \cite{Gu}.

\

\noindent{\bf Theorem B.} \cite{GeWang_Proc} {\it Let $(M^4, g)$ is a $4$-dimensional closed Riemannian manifold with non-negative
scalar curvature, then
\begin{equation}\label{eq1_add1}
\int_M |Ric-\frac{\overline R } 4 g|^2 dv (g)\le 4\int_M |Ric-\frac{ R } 4 g|^2dv(g),
\end{equation}
where $\overline R=vol (g)^{-1} \int_M R dv(g)$ is the average of the scalar curvature $R$ of $g$. Or
equivalently, we have
\begin{equation}\label{eq3_add}
\frac 83 vol(g) \int_M\s_2(g) dv(g)\le \left(\int_M \s_1(g) dv(g)\right)^2.
\end{equation}
}
\

In fact, one can find inequality (\ref{eq3_add}) in the argument of Gursky \cite{Gu}.
This argument uses a crucial property of $\s_2$-scalar curvature 
 that
$\int_{M}\s_2(g) dv(g)$ is a conformal invariant, which is only true on $4$-dimensional manifolds. Nevertheless, inspired  by our previous work in \cite{GLW} we conjectured in \cite{GeWang_Proc} that this is true for $3$-dimensional
manifolds.
In this paper, by improving the analysis developed in the study of $\s_k$-Yamabe problem,  we
give an affirmative answer to this conjecture. Namely we will show that  Theorem A holds under the condition of non-negativity of the scalar
curvature for dimension $n=3$. 

\begin{thm} \label{mainthm1} Let $(M^3, g)$ is a $3$-dimensional closed Riemannian manifold with non-negative
scalar curvature. We have 
\begin{equation}\label{eq_1}
\int_M |Ric-\frac{\overline R } 3 g|^2 dv (g)\le 9 \int_M |Ric-\frac{ R } 3 g|^2dv(g).
\end{equation}
Moreover, equality holds if and only if $(M^3, g)$ is a space form.
\end{thm}

Without the condition of non-negativity of  the scalar curvature, Theorem \ref{mainthm1} is not true. Examples can be found in
\cite{DT}. When $n>4$, Theorem A is also not true under a weaker condition that the scalar curvature is positive.
For various problems related to the De Lellis-Topping inequality, see \cite{GWX}.

Our proof is based on the study of a new conformal invariant.
From now, let $n=3$. We define
\begin{equation}
\ds \widetilde Y([g_0]):=\sup_{ g\in {\cal C}_1([g_0])}\frac {\ds vol(g)\int_M \s_2( g) dv( g)}
{\ds (\int_M \s_1( g) dv( g))^2},
\end{equation}
where ${\cal C}_1([g_0]):=\{g=e^{-2u}g_0\,|\, R>0\}$ and $[g_0]:=\{g=e^{-2u}g_0\}$.
We define  the first Yamabe constant  on 3-dimensional manifolds by
\[Y_1([g_0]):= \inf_{\tilde g \in [g_0]} \frac {\ds \int_M\s_1(\tilde g) dv(\tilde g)}
{\ds (vol(\tilde g))^{\frac {1}{3}}}.\]
Since $\s_1(g)=R/2(n-1)$, the first Yamabe constant $Y_1([g_0])$ is a positive constant multiple of the ordinary Yamabe constant.
Theorem \ref{mainthm1} follows from 
the observation mentioned above and the following
\begin{thm} \label{mainthm2}  Let $(M^3, g)$ is a $3$-dimensional  closed Riemannian manifold with positive Yamabe constant  $Y_1([g_0])>0$, then
\begin{equation}\label{Invariant}
\ds \widetilde Y([g_0])\le \frac 1 3.
\end{equation}
\end{thm}

To show Theorem \ref{mainthm2} we will study a fully nonlinear
Yamabe type equation \eqref{new}, which is closely related to the
$\s_k$-Yamabe problem initiated  in \cite{V}, \cite{CGY} and studied
by many mathematicians. (See for example \cite{guan} and
\cite{Survey}) Though the fully nonlinearity, the corresponding $\s_k$ Yamabe
equation  shares very nice properties. (See \cite{GuanWang1} and
\cite{LL}) A nice application of the analysis developed in the study
of the $\s_k$-Yamabe problem is the 4-dimensional sphere theorem
obtained by Chang-Gursky-Yang in \cite{CGY2}. As another
application, with C.-S. Lin we obtained in \cite{GLW} a
3-dimensional sphere theorem. Another proof was given by
Catino-Djadli  in \cite{CD}. See also \cite{CGY}, \cite{GeWang2}, \cite{GuanLinWang}, \cite{GuV}, \cite{GuV3},  \cite{Wang},
and especially a survey paper \cite{Gu2}
for other applications.
This paper can be seen as a new
application of this analysis. However, comparing to the ordinary Yamabe
problem and $\s_k$-Yamabe problem we encounter an extra difficulty,
without a corresponding Sobolev inequality, which is in fact
inequality \eqref{Invariant} that we want to prove.

Theorem \ref{mainthm1} gives also
a new characterization of three-dimensional spherical space forms. Another related characterization of three-dimensional  space forms
was recently given by Gursky and Viaclovsky in \cite{GuV2}.

The paper is organized as follows. In Section 2 we consider the new
conformal invariant $\widetilde Y$ and its related energy
functional. The critical point of this energy functional satisfies a
Yamabe type equation \eqref{new} below. We show in Lemma \ref{lem2}
that any critical point satisfies $\widetilde Y\le 1/3$. Hence  to
prove Theorem \ref{mainthm2} we only need  to prove that  $\widetilde Y$
is achieved. This is in fact a new Yamabe type problem, with a new
difficulty -without a corresponding Sobolev inequality. This problem
is difficult and still remains open. Instead of  attacking this
problem directly we  consider a suitable perturbed problem. 
 This perturbed equation, to find it is a very delicate issue,
 is introduced in Section 3. 
In Section 4 we prove  first local $C^2$ estimates and then global $C^2$ estimates for the flow, by using the local estimates.
The uniform parabolicity of the flow is proved in Section 5. One of key estimates (Lemma \ref{lastlem}) and main Theorems are proved in Section 6.
Related problems and Conjectures are proposed in Section 7.

\section{A new conformal invariant and a related flow}

\

Let us first recall the definition of the $k$-scalar curvature, which was first introduced by Viaclovsky \cite{V}
and has been intensively studied by many mathematicians, see for example the references in \cite{GLW} and two survey papers \cite{guan} and \cite{Survey}. 
Let
\[S_g=\frac 1{n-2}\left(Ric-\frac {R}{2(n-1)}\cdot g\right)\]
be the Schouten tensor of $g$. For an integer $k$ with $1\le k\le n$ let $\s_k$ be the $k$-th elementary symmetric function in $\R^n$.
The $k$-scalar curvature is defined by
\[\s_k(g):=\s_k(\L_g),\]
where $\L_g$ is the set of eigenvalue of the matrix $g^{-1}\cdot
S_g$. In particular,
\[\begin{array}{rcl}
\s_1(g)=\ds \frac {R}{2(n-1)}, \quad \quad 
\s_2(g)&=& \ds \frac 1{2(n-2)^2} \left\{-|Ric|^2+\frac n{4(n-1)} R^2\right\}.\\
\end{array}\]
We have in \cite{GeWang_Proc} the following observation.
\begin{lem}\label{keylemma} (\cite{GeWang_Proc}) Inequality (\ref{eq1}) is equivalent to
\[\left(\int_M \s_1(g) dv(g)\right)^2\ge \frac {2n}{n-1} vol(g) \int_M\s_2(g) dv(g).\]
\end{lem}

\


Let $g_0$ be a metric on $M^3$ with positive scalar curvature and
${\cal C}_1([g_0]):=\{ g\in [g_0]\,|\, \s_1(g)>0\}$. Define an
energy functional
\begin{equation}
\ds {\cal E}(g):=\frac {\ds vol(g)\int_M \s_2( g) dv( g)} {\ds
(\int_M \s_1( g) dv( g))^2}.
\end{equation}
and
\begin{equation}\ds \widetilde Y([g_0]):=\sup_{ g\in {\cal C}_1([g_0])}{\cal E}(g).
\end{equation}
$\widetilde Y([g_0])$ is a new conformal invariant. To show Theorem \ref{mainthm2} is equivalent to show that
this invariant is always less than or equal to $1/3$.
A critical point of ${\cal E}$ in   $ {\cal C}_1([g_0]) $ satisfies a new Yamabe type equation

\begin{equation}\label{new}
\frac {\s_2(g)-
3r_2(g)}{\s_1(g)}=-2s(g),
\end{equation}
where $r_2(g)$ is the average of $\s_2(g)$  and $s(g)$ the average of
$\frac{\s_2(g)}{\s_1(g)}$ with respect to the measure $\s_1(g)dv(g)$
are defined by
\[r_2(g):=\frac {\ds\int_M \s_2(g)dv(g)}{vol(g)} \quad
\hbox{ and } \quad
s(g):=\frac {\ds\int_M \s_2(g) dv(g)}{\ds \int_M\s_1(g) dv(g)}.
\]
We observe that solutions of \eqref{new} have an interesting property.
\begin{lem}\label{lem2}
 Every solution $g\in {\cal C}_1([g_0])$ of \eqref{new} satisfies
\begin{equation}\label{13}
{\cal E}(g)\le 1/3,\end{equation}
and equality if and only if $g$ is an Einstein metric.
\end{lem}
\pr From the Newton inequality $$\ds\frac{\s_2(g)}{\s_1(g)}\le \frac13 \s_1(g),$$
we have
$$
\begin{array}{lll}
\ds\frac13 (\int_M \s_1(g)dv(g))^2&\ge \ds\int_M
\frac{\s_2(g)}{\s_1(g)}dv(g)\int_M \s_1(g)dv(g)\\
&=\ds 3r_2(g)\int_M \frac{1}{\s_1(g)}dv(g)\int_M
\s_1(g)dv(g)-2s(g)vol(g)\int_M \s_1(g)dv(g)\\
&\ds\ge 3r_2(g)(vol(g))^2-2s(g)vol(g)\int_M \s_1(g)dv(g)\\
&\ds=\int_M \s_2(g)dv(g) vol(g).
\end{array}
$$
In the first equality we have used Equation \eqref{new} and in the second inequality
the Cauchy-Schwarz inequality.
It is clear to see that equality holds if and only if
$$\ds\frac{\s_2(g)}{\s_1(g)}= \frac13 \s_1(g),$$
and hence if and only if $(M^3,g)$ is an Einstein manifold. \qed

\

Therefore, to prove Theorem \ref{mainthm2} we only need to prove the existence of the maximum
of functional $\cal E$ in ${\cal C}_1$. This is a new Yamabe type problem.
However to prove the existence of the the maximum
of functional $\cal E$ is very difficult. One would meet  not only the typical difficulty -loss of the compactness-
of the ordinary Yamabe problem (and many other geometric variational problems, for example, harmonic maps, Yang-Mills fields), the fully nonlinearity
of the $\s_k$-Yamabe problem, but also a new problem that we have not  a corresponding
(optimal) Sobolev inequality yet. This corresponding Sobolev inequality is
\[\sup {\cal E}(g) < \infty, \hbox{ or } \sup {\cal E}(g) \le \frac 13.\]
This is in fact what we want to show.
Hence we  need to consider certain suitable perturbed functionals.

\section{A perturbed problem and its flow}
As mentioned in the Introduction, to find a suitable perturbed problem  is a delicate issue. 
 Let
$\varepsilon>0$ be some small constant and $ g\in {\cal C}_1([g_0])$.
We define
\begin{equation}
\ds {\cal E}_\e(g):=\frac {\ds\left(\int_M (\s_2( g)-\e e^{4u}) dv(
g)\right)\left(\int_M  e^{\e u} dv( g)-\e (\int_M \s_1( g) dv(
g))^{3-\e}\right)  } {\ds \left(\int_M (\s_2( g)-\frac{\e}{2}
e^{4u}) dv( g)\right)^\e\left(\int_M \s_1( g) dv( g)\right)^2}.
\end{equation}
This perturbed function is well defined in a smaller space 
$$
\begin{array}{ll}
 {\cal C}_{1,\e}([g_0])
 :=\\
 \ds\{g\in  {\cal C}_1([g_0])|\quad \int_M (\s_2( g)-\e e^{4u}) dv(
 g)>0,
 \quad\int_M  e^{\e u} dv( g)-\e (\int_M \s_1( g) dv( g))^{3-\e}>0\}.
 \end{array}
$$
This functional looks quite complicated. But it satisfies all properties we want to have.
Denote the maximum of ${\cal E}_\e$ in ${\cal C}_{1,\e}([g_0])$ by
$$
M_\e:=\sup_{g\in {\cal C}_{1,\e}([g_0])}{\cal E}_\e(g).
$$
For this perturbed energy functional, the corresponding
Euler-Lagrange equation could be written as follows
\begin{equation}\label{p-eq}
 \frac {\s_2(g)- (\nu_1(g)e^{\e
u}+\nu_2(g)e^{4 u})}{\s_1(g)}+\mu(g)=0
\end{equation}
where $\nu_1(g)$, $\nu_2(g)$ and $\mu(g)$ are given respectively
$$
\nu_1(g):=(3-\e)\frac{k(g)\ds(\int_M (\s_2( g)-\e e^{4u}) dv(
g))}{\ds\int_M  e^{\e u} dv( g)-\e (\int_M \s_1( g) dv( g))^{3-\e}},$$
$$
\nu_2(g):=\e k(g)\left(1-\frac{\ds\e \int_M (\s_2( g)-\e e^{4u}) dv(
g)}{\ds2\int_M (\s_2( g)-\frac{\e}{2} e^{4u}) dv( g)}\right),
$$
$$
\mu(g):=k(g)\ds \int_M (\s_2( g)-\e e^{4u}) dv(g)
\left(\frac{2}
{\ds\int_M \s_1( g) dv( g)}+\frac{\e (3-\e)(\int_M \s_1( g) dv( g))^{2-\e}}{\ds\int_M
e^{\e u} dv( g)-\e (\int_M \s_1( g) dv( g))^{3-\e}}\right),
$$
with
$$\begin{array}{rcl}
k(g)&:=&\ds\vs \frac{\ds\int_M (\s_2( g)-\frac{\e}{2} e^{4u}) dv(
g)}{\ds\int_M (\s_2( g)-\frac{\e}{2} e^{4u}) dv( g)-\e\int_M (\s_2(
g)-\e e^{4u}) dv(
g)}\\
&=&\ds\vs \frac{\ds\int_M (\s_2( g)-\frac{\e}{2} e^{4u}) dv( g)}
{(1-\e)\ds\int_M(\s_2( g)-\frac{\e(1-2\e)}{2(1-\e)} e^{4u}) dv(
g)}\ge 1.
\end{array}
$$
By definition we have
\begin{lem}\label{lem3} We have
\begin{itemize}
\item[(i)] $\frac{\nu_1(g)}{\mu(g)}\le\frac1{\e(\int\s_1(g)dv(g))^{2-\e}}
$.
\item[(ii)]$\widetilde Y([g_0])\le \limsup_{\e\to 0} M_\e$.
\end{itemize}
\end{lem}
\pr The proof is easy to check.\qed

We want to show that
\begin{itemize}
 \item[(1)] $M_\e$ is achieved by some $g_\e \in {\cal C}_{1,\e}([g_0])$ for $\e>0$, which  certainly satisfies
\eqref{p-eq}.
\item[(2)] Every solution $g$ of \eqref{p-eq} satisfies an estimate
\begin{equation}\label{est-add}
{\cal E}_\e(g) \le (\frac 2{C\e})^\e \frac 1{3(1-\e)},
               \end{equation}
  where $C$ is a constant independent of $\e$. 
\end{itemize}
This implies that $\widetilde Y([g_0])\le \limsup_{\e\to 0} M_\e\le 1/3$.
Estimate \eqref{est-add} will be proved in Lemma \ref{lastlem} below.
To study the achievement of $M_\e$, we introduce a conformal flow, which
is different from the Yamabe flow considered in \cite{GLW}.
\begin{equation} \label{flow}
\frac {du}{dt}=-\frac 12g^{-1}\cdot \frac d{dt}g:=e^{-2u}\frac {\s_2(g)-
(\nu_1(g)e^{\e u}+\nu_2(g)e^{4 u})}{\s_1(g)}+\mu(g)e^{-2u}+m(g),
\end{equation}
where $m(g)$ is chosen by
$$
\int_M \s_1(g)\left( e^{-2u}\frac {\s_2(g)-
(\nu_1(g)e^{\e u}+\nu_2(g)e^{4 u})}{\s_1(g)}+\mu(g)e^{-2u}+m(g)  \right)dv(g)=0.
$$

\begin{pro}\label{pro} Let $n=3$. Flow (\ref{flow}) preserves $\ds\int_M \s_1(g) dv(g)$, while it
increases ${\cal E}_\e(g)$, provided $g(t)\in  {\cal C}_{1,\e}([g_0])$.
\end{pro}

\noindent{\it Proof.}
It is clear that the flow preserves $\int_M\s_1(g) dv(g)$. By a direct computation we have

\[
\begin{array}{rcl}
\ds\vs \frac d{dt} {\cal E}_\e(g) &=&\ds  \frac{\ds{\cal
E}_\e(g)\int_M e^{-2u}\s_1(g) \left(\frac {\s_2(g)- (\nu_1(g)e^{\e
u}+\nu_2(g)e^{4 u})}{\s_1(g)}+\mu(g)\right)^2}{\ds k(g)\int_M (\s_2(
g)-\e e^{4u}) dv( g)} \ge 0.
\end{array}
\]
\qed

Since the flow increases ${\cal E}_\e(g)$, the flow preserves the properties
$\int_M (\s_2( g)-\e e^{4u}) dv(
 g)>0$, $ \int_M  e^{\e u} dv( g)-\e (\int_M \s_1( g) dv( g))^{3-\e}>0$.
We will show below that the flow preserves ${\cal C}_1([g_0])$, and hence ${\cal C}_{1,
\e}([g_0])$.
  This is certainly one of crucial properties
of the flow.


\section{$C^2$ estimates}
In this  section, we will establish a priori estimates for flow
(\ref{flow}). Local estimates for this class of fully nonlinear
conformal equations were first given in \cite{GuanWang1}. Since then there
are many extensions. See for instance \cite{Chen} and the survey
paper \cite{Survey}. Let $ \Gamma_k^+$ be a convex open cone -the Garding cone- defined by
\[\Gamma_k^+=\{\Lambda=(\l_1,\l_2,\cdots, \l_n)\in \R^n\,|\,
\sigma_j(\Lambda)>0, \forall j\le k\}.\] 
Similarly, we say a symmetric matrix $W\in \Gamma_k^+$ if the set of eigenvalues of $W$ belongs to $\Gamma_k^+$. 
By $g\in {\cal C}_k^+$ we
mean that 
$g^{-1}\cdot S_g (x)$ belongs to $\Gamma_k^+$ for any $x\in
M$.  
If $g=e^{-2u} g_0$, we have the transformation formula of the Schouten tensor
\[ S_g=\uuu.\]
Therefore, $g=e^{-2u} g_0\in {\cal C}_k$ if and only if
\[ (\uuu)(x)\in \Gamma_k^+, \quad \forall x\in M.\]

To establish  a priori estimates, we first need a technical key lemma.\\

\begin{lem}
\label{lem3.1} For $1<k\le n$ set $F=\frac{\s_{k}}{\s_{k-1}}$. We have
\begin{itemize}
\item[1)]  the matrix  $(F^{ij})(W)$ is semi-positive definite at $W\in \Gamma_{k-1}^+$
 and is   positive definite at $W \in \Gamma_{k-1}^+\backslash {\cal
 R}_1$, where ${\cal
 R}_1$  is the set of  matrices of rank $1$.\\
\item [2)] The function $F$ is concave in the cone $\Gamma_{k-1}^+$.
When $k=2$,  for all $W\in \Gamma_1^+$ and for all $R=(r_{ij})\in
{\cal S}_n$, we have \beq\label{eq3.1} \ba{rcl} \ds\vs
\sum_{ijkl}\frac {\p^2}{\p w_{ij}\p w_{kl}}\left(\frac
{\s_2(W)}{\s_1(W)}\right)r_{ij}r_{kl} &=&\ds -\frac{
\sum_{ij}(\s_1(W)r_{ij}- \s_1(R)w_{ij})^2 }{ \s_1^3(W)}. \ea \eeq
\end{itemize}
\end{lem}
\pr For the proof, see \cite{GLW}. \qed

\

Assume $g_1\in  {\cal C}_{1,\e}([g_0])$.  We consider flow \eqref{flow}
with the initial metric $g_1$.
 Lemma \ref{lem3.1} implies that 
(\ref{flow}) is  parabolic. By the standard implicit function
theorem we have the short-time existence result. Let $T^*\in
(0,\infty]$ so that $[0, T^*)$ is the maximum interval for
 the existence of the flow $g(t)\in {\cal C}_{1,\e}([g_0])$.

\begin{thm}\label{thmlocal1} Assume that $n= 3$,  and  $g(0)=g_1\in {\cal C}_{1,\e}([g_0])$. Let $u$ be a solution of
(\ref{flow}) in a geodesic ball $B_R\times[0, T]$ for $T<T^*$ and
$R<\tau_0$, the injectivity radius of $M$.  Then there is a constant
$C$ depending only on $(B_R, g_0)$  and independent of $T$ such that
for any $(x,t) \in B_{R/2}\times[0, T]$
\begin{equation}\label{eq_thm1.1}
\ba{rcl}
|\n u|^2+|\n^2 u|&\le&  \ds\vs
C(1+\frac{\nu_1(g)}{\mu(g)}e^{-(2-\e)\inf_{B_R}u})\\
&\le& \ds C(1+\frac1{\e(\int\s_1(g)dv(g))^{2-\e}} e^{-(2-\e)\inf_{B_R}u}).
\ea
\end{equation}
\end{thm}

\pr  In the proof, $C$ (resp. $c$) is a constant independent of $T$, which may
vary from line to line. Let $W=(w_{ij})$ be an $n\times n$ matrix
 with $$w_{ij}=u_{ij} +u_iu_j-\frac {|\n u|^2}2
 (g_0)_{ij}+(S_{g_0})_{ij}.$$ Here $u_i$ and $u_{ij}$ are the first
 and second derivatives of $u$ with respect to the background metric $g_0$.
 Define
 $$\nu:=\nu_1(g)e^{-(4-\e)u}+\nu_2(g), \quad \bar \nu:=\nu_1(g)e^{-(4-\e)u}$$ and
 $$F(W,u):=\ds\frac{\s_2(W)-\nu}{\s_1(W)}.$$
Set
 \beq \ba{rcl}\ds\vs
(F^{ij}(W,u))&:=&\ds \left(\frac{\partial F}{\partial w_{ij}}(W)\right)\\
&=&\ds\vs \left(\frac{\s_1(W)T^{ij}-\s_2(W)\delta^{ij} +\nu
\delta^{ij}}{\s_1^2(W)}\right) \ea \eeq where $(T^{ij})=(\s_1(W)
\delta^{ij} -w^{ij})$ is the first Newton transformation associated
with $W$, and $\delta^{ij}$ is the Kronecker symbol. From Proposition \ref{pro}, we know $\nu_1(g)>0$, $\nu_2(g)>0$ and
$\mu(g)>0$. In view of
Lemma \ref{lem3.1} we know that $(F^{ij})$ is positive definite and
$F$ is concave in $W\in \Gamma_1^+$. Moreover, we have
\begin{equation}
\label{concave}
\sum_{ijkl}\frac {\p^2\left( F(W,u)\right)}{\p w_{ij}\p w_{kl}} r_{ij}r_{kl} \le -2\frac{\nu (\sum_{i}r_{ii})^2}{ \s_1^3(W)}.
\end{equation}
Let $S(TM)$ denote the unit tangent
bundle of $M$ with respect to the background metric $g_0$. We define
a function $\tilde G:\;S(TM)\times [0, T]\;\to\R$ \beq \tilde
G(e,t)= (\n^2 u+|\n u|^2 g_0)(e,e).\eeq 
Without loss of generality,
we assume $R=1$. Let $ \rho\in C^\infty_0(B_1)$ be a cut-off
function defined as in \cite{GuanWang1} such that
\begin{equation}\begin{array}{rcll}\label{eq4.1}
\vs \ds \rho & \ge & 0, & \hbox{ in } B_1,\\
\vs\ds \rho & =& 1, &\hbox{ in } B_{1/2},\\
\vs\ds |\n \rho (x)|& \le &2 b_0 \rho^{1/2}(x), & \hbox{ in } B_1,\\
|\n^2 \rho| & \le  & b_0,   & \hbox{ in } B_1.
\end{array}\end{equation}
Here $b_0>1$  is a constant. Since $e^{-2u}g_0\in  {\cal C}_1$, to
bound $|\n u|$ and $|\n^2 u|$ we only need to bound
 $(\n^2 u+|\n u|^2 g_0)(e,e)$ from above for all $e\in S(TM) $ and for all $t\in [0, T] $. For this purpose, 
consider $G(e,t)=\rho(x)\tilde G(e,t) $.
 Assume $(e_1,t_0)\in S(T_{x_0}M)\times (0, T]$ such that
\beqr \label{eq4.2}
\vs \ds G(e_1,t_0)=\max_{S(TM)\times [0, T]} G(e,t).\eeqr
We may further assume that \beqr \label{eq4.4}\ds G(e_1,t_0)>n
\max_{B_1} \s_1(g_0). \eeqr Let $(e_1,\cdots,e_n)$ be an orthonormal
basis at point $(x_0,t_0)$.
Now choose the normal
coordinates around $x_0$ such that at point $x_0$
\[
\frac{\p}{\p x_1}=e_1
\]
and consider the function $G$ on $M\times [0,T]$ defined by
\[
G(x,t):=\rho(x)(u_{11}+ |\n u|^2)(x,t).
\]
Clearly, $(x_0,t_0)$ is a maximum point of $G(x,t)$ on $M\times
[0,T]$.
At $(x_0, t_0)$, we have
\begin{eqnarray}
\label{y1}
0&\le & G_t=\rho(u_{11t}+2\sum_{l}u_lu_{lt}),\\
\label{y2} 0&=&G_j=\frac{ \rho_j }{\rho} G + \rho
 (u_{11j}+2\sum_{l\ge 1}u_l u_{lj}),
\quad \hbox{ for any } j,\\
0&\ge &(G_{ij} )=\ds \left(\frac{ \rho
\rho_{ij}-2\rho_i\rho_j}{\rho^2} G+ \rho  (u_{11ij}+\sum_{l\ge
1}(2u_{li}u_{lj}+2 u_l u_{lij}))\right).
\end{eqnarray}
Recall that  $(F^{ij})$ is  definite positive. Hence, we have \beq
\label{x-2} \ba{rcl} 0& \ge &\ds\vs  \sum _{i,j\ge 1}F^{ij}G_{ij}
-G_t\\
& \ge &  \ds\vs\sum_{i,j\ge 1}
 F^{ij} \frac{ \rho \rho_{ij}-2\rho_i\rho_j}{\rho^2} G +
\rho \sum_{i,j\ge 1} F^{ij}(u_{11ij}+\sum_{l\ge 1}(2u_{li}u_{lj}+2 u_l u_{lij}))\\
&& \ds-\rho(u_{11t}+2\sum_{l\ge 1}u_lu_{lt}). \ea \eeq First, from
the definition of
 $\rho$, we have
\beq \label{x-1} \sum_{i,j\ge 1}
 F^{ij} \frac{ \rho \rho_{ij}-2\rho_i\rho_j}{\rho^2} G \ge -C
 \sum_{i,j \ge 1}|F^{ij}|\frac 1 \rho  G,
\eeq and \beq \label{x-3} \ba{rcl}
\ds\sum_{i,j\ge 1} |F^{ij}|  & \ge &\ds\vs\sum_{i} F^{ii}  \\
&=& \ds \vs \left( n-1-\frac {n\s_2(W)}{\s_1^2(W)}\right)+ \frac
{n\nu}{\s_1^2(W)}\ge C \sum_{i,j\ge 1} |F^{ij}|, \ea \eeq since $W$
is positive definite. From \eqref{x-3} we have \beq\label{x-4}
\sum_i F^{ii} \ge \frac {n-1}
2+\frac{n\nu}{\s_1^2(W)}=1+\frac{3\nu}{\s_1^2(W)}. \eeq Using the
facts that
\begin{equation}
\label{0x1}
 u_{kij}=u_{ijk}+\sum_m R_{mikj}u_m,
\end{equation}
\begin{equation}
\label{0x2}
 u_{kkij}=u_{ijkk}+
 \sum_m (2R_{mikj}u_{mk}-Ric_{mj}u_{mi}-Ric_{mi}u_{mj}
 -Ric_{mi,j}u_m+R_{mikj,k}u_m)
\end{equation}
and
\begin{equation}
\label{0x3} (\sum_l u_l^2)_{11} =2\sum_l(u_{11l}u_l+ u_{1l}^2) +O
(|\n u|^2),
\end{equation}
we have \beq\label{x1} \ba{rcl}
  \ds \sum_{i,j\ge 1} F^{ij} u_{11ij}
&\ge & \ds\vs \sum_{i,j\ge 1} F^{ij} \left(w_{ij11}-(u_{11})_iu_j -u_i (u_{11})_j+\sum_{l\ge 1} (u_{1l}^2 + u_{11l}u_l)(g_0)_{ij}\right)\\
 & & \ds\vs -2\sum_{i,j\ge 1} F^{ij}u_{i1}u_{j1}-
 C(1+|\n^2 u|+|\n u|^2)\sum_{i,j\ge 1} |F^{ij}|
\ea\eeq and \beq\label{x2}\ba{rcl} \ds\vs \sum_{i,j, l} F^{ij}
u_{l}u_{lij}
 &\ge & \ds \sum_{i,j, l} F^{ij} u_{l}w_{ijl}-\sum_{i,j, l} F^{ij}(u_lu_{il}u_j+u_lu_iu_{jl}) \\
 &&\ds +\frac1 2\sum_{i,j} F^{ij}\langle \n u, \n(|\n u|^2)\rangle (g_0)_{ij} -C(1+|\n u|^2)\sum_{i,j\ge 1} |F^{ij}|.
\ea\eeq Combining (\ref{x1}) and (\ref{x2}), we deduce
\beq\label{x3}\ba{rcl} & &\ds\vs \sum_{i,j\ge 1}
F^{ij}(u_{11ij}+2\sum_{l\ge 1}(u_{li}u_{lj}
+ u_l u_{lij}))\\
&\ge &\ds\vs \sum_{i,j\ge 1} F^{ij} (w_{ij11}+2\sum_{l\ge 1} w_{ijl}u_l) +2\sum_{i,j\ge 1} F^{ij} \sum_{l\ge 2}u_{li} u_{lj}+\sum_{i,j,l\ge 1} u_{1l}^2 F^{ij}(g_0)_{ij}\\
& &\ds\vs -\sum_{i,j} F^{ij}\left[(u_{11}+ |\n u|^2)_iu_j+ u_i(u_{11}+ |\n u|^2)_j-\langle \n u, \n(u_{11}+ |\n u|^2)\rangle (g_0)_{ij}\right]\\
&&\ds\vs -C(1+|\n^2 u|+|\n u|^2)\sum_{i,j\ge 1} |F^{ij}|\\
&\ge&\ds \vs \sum_{i,j} F^{ij} (w_{ij11}+2\sum_l w_{ijl}u_l)+u_{11}^2\sum_{i,j} F^{ij}(g_0)_{ij}\\
&&\ds \vs + \sum_{i,j} F^{ij}\left( \rho_iu_j+\rho_ju_i
-\langle\n\rho,\n u\rangle (g_0)_{ij} \right)
\frac{G}{\rho^2}-C(1+|\n^2 u|+|\n u|^2)\sum_{i,j\ge 1} |F^{ij}|.
\ea\eeq
In the last inequality we have used \eqref{y2}.
 Now, we want to estimate $ \sum_{i,j,l} F^{ij} w_{ijl}u_l$
and $\sum_{i,j} F^{ij} w_{ij11}$ respectively.
By differentiating $F$ we get
 \beq\label{eq4.11}\ba{llll}
\ds \sum_l F_lu_l= \sum_{i,j,l}F^{ij} w_{ijl}u_l+\sum_{l} \frac{\p F}{\p u
}u_l^2&=\ds\sum_{i,j,l} F^{ij} w_{ijl}u_l+\sum_{l} \frac{(4-\e)\bar\nu
u_l^2}{\s_1(W) }. \ea\eeq
By differentiating $F$ twice and
using the concavity  (\ref{concave}) of $F$ in $W$, we have
\beq\label{eq4.12x}\ba{rcl} \ds \sum_{i,j} F^{ij} w_{ij11}&=& \ds
F_{11}-\sum_{i,j,k,m}\frac{\p^2 F}{\p w_{ij}\p w_{km}}
w_{ij1}w_{km1}\\
&&-\ds 2\sum_{i,j}\frac{\p^2 F}{\p w_{ij}\p u}
w_{ij1}u_1-\frac{\p^2 F}{\p^2 u}
u_{1}^2-\frac{\p F}{\p u}
u_{11}\\
&\ge&\ds F_{11}+\frac{2\nu (\sum_i
w_{ii1})^2}{(\s_1(W))^3}+\frac{2(4-\e)\bar \nu (\sum_i
w_{ii1})u_1}{(\s_1(W))^2} +
\frac{(4-\e)^2\bar\nu u_1^2}{\s_1(W)}- \frac{(4-\e)\bar\nu u_{11}}{\s_1(W)}\\
&\ge&\ds F_{11}+ \frac{(4-\e)^2\bar\nu u_1^2}{2\s_1(W)}-
\frac{(4-\e)\bar\nu u_{11}}{\s_1(W)}. \ea \eeq These
estimates give \beq\label{eq4.12}\ba{rcl} \ds \sum_{i,j\ge 1}
 F^{ij}
(w_{ij11}+2\sum_{l\ge 1} w_{ijl}u_l)&\ge& F_{11} +2 \sum_l F_lu_l\\
&&\ds -\frac{(4-\e)\bar\nu u_{11}}{\s_1(W)}-\sum_{l=1}^n
\frac{2(4-\e)\bar \nu u_l^2}{\s_1(W) }. \ea\eeq Recall from
(\ref{flow}) that \beq\label{eq4.13} F=u_t-\mu(g) e^{-2u}-m(g). \eeq
Hence we have
\beq\label{eq4.14} F_{11}=u_{11t}-\mu(g)
e^{-2u}(-2u_{11}+4u_1^2), \eeq \beq\label{eq4.15}
F_{l}=u_{lt}-\mu(g) e^{-2u}(-2u_{l}),\;\forall\, l=1,\cdots,n. \eeq
Gathering 
(\ref{x-2}), (\ref{x-1}), (\ref{x-3}), (\ref{x3}) (\ref{eq4.12}),
(\ref{eq4.14}) and (\ref{eq4.15}), we obtain \beq \label{eq4.16}
\ba{rcl} 0 &\geq& \ds\vs -C\left(\sum_{i,j}
\left|F^{ij}\right|\right) \frac{G}{\rho}+\rho\left(\sum_i F^{ii}
\right)u_{11}^2 -C
\rho\left(\sum_{i,j} \left|F^{ij}\right|\right) (1+|\n u|^2+|\n^2 u|)\\
&& \ds\vs + \sum_{i,j} F^{ij}\left( \rho_iu_j+\rho_ju_i
-\langle\n\rho,\n u\rangle (g_0)_{ij} \right)
\frac{G}{\rho} -\rho \mu(g) e^{-2u}\left(-2u_{11}-4\sum_{l=2}^n u_l^2\right)\\
&&\ds -\rho \frac{(4-\e)\bar \nu }{\s_1(W)}(u_{11}+ 2\sum_{l=1}^n
u_l^2). \ea \eeq 
 From the fact $W\in \Gamma_1^+$, we have that $u_{11}(x_0,t_0)\ge\frac{1}{20}|\nabla
u|^2(x_0,t_0)$, and hence  $G(x_0,t_0)\le 21\rho(x_0) u_{11}(x_0,t_0)$
(see (44) in \cite{GLW}). Multiplying \eqref{eq4.16} by $\rho$ we deduce
\beq \label{estimate1}0\ge \sum_{i} F^{ii}(-C G +
(\frac{G}{21})^2-CG^{\frac 32}) +\rho e^{-2u}(\mu (g)\frac{2G}{21}- 8
\nu_1(g) e^{-(2-\e)u}\frac{G}{\s_1(W)}). \eeq When
$\frac{G}{\s_1(W)}\ge 2352=16\times (21)^2/3$, it follows from
(\ref{x-4}) that
\[ \frac12 \sum_{i} F^{ii}(\frac{G}{21})^2 -8 \nu_1(g)\rho
 e^{-(4-\e)u}\frac{G}{\s_1(W)}\ge \nu
(\frac{G^2}{294\s_1^2(W)} -8 \rho \frac{G}{\s_1(W)})\ge 0.\]
Together with \eqref{estimate1}, we have  $$0\ge \sum_{i} F^{ii}(-C
G + \frac12(\frac{G}{21})^2-CG^{\frac 32}), $$ from which we easily have $$
G(x_0,t_0)\le C.
$$
This gives the desired result.
When $\frac{G}{\s_1(W)}< 2352$, the desired result follows from \eqref{estimate1} and Lemma \ref{lem3} (i). 
 \qed

\begin{rem}\label{thmlocal2} Let $g=e^{-2u}g_0\in {\cal C}_{1,\e}$ be a solution of
(\ref{p-eq}) in a geodesic ball $B_R$ and
$R<\tau_0$, the injectivity radius of $M$. Then there is a  constant
$C$ depending only on $(B_R, g_0)$  such that for any
 $x \in B_{R/2}$ the estimate (\ref{eq_thm1.1}) holds.
\end{rem}

\begin{cor}\label{thmglobal1}
Under the same assumptions as in Theorem 
\ref{thmlocal1},
 there is a  constant $C$ depending only on  $g_0$  (independent of $T$) such that
for any  $t\in [0, T]$
\begin{equation}
\|u\|_{C^2(M)} \le C.
\end{equation}
\end{cor}

\pr By Proposition \ref{pro}, we may assume that  $\ds \int_M \s_1( g)
dv( g)\equiv 1$ without loss of generality. Thus, we have a uniform
volume bound, namely \beq\label{volumebound} vol(g)\le
(Y_{1}([g_0]))^{-3}. \eeq


\noindent{\bf Claim.} There is a constant $C>0$ independent of $T\in
[0, T^*)$ such that $\forall t\in [0,T]$
\begin{equation}\label{eq6.7}
u(t,x)\ge C.
\end{equation}
Set $m(t)=\min_{x\in M} u(t,x)$ and $u(t,x_t)=m(t)$. 
We prove the claim by a
contradiction argument and assume that there exists a sequence
$\{t_n\}$ such that $t_n\to T$ and $m(t_n)\to -\infty$. Applying
Theorem \ref{thmlocal1}, we have for all $x\in M$ and $n\in\N$
$$
|\nabla u (t_n,x)|^2\le \frac{C}{\e} e^{-(2-\e)m(t_n)},
$$
which implies for all $x\in B(x_{t_n}, \sqrt{\e}
e^{(1-\e/2)m(t_n)})$
$$
|u(t_n,x)-m(t_n)|\le C.
$$
As a consequence, we infer
$$
vol(g(t_n))\ge C \int_{B(x_{t_n}, \sqrt{\e} e^{(1-\e/2)m(t_n)})}
e^{-3m(t_n)} dv(g_0)\ge\e^{3/2} e^{(-3\e/2)m(t_n)}\to\infty
$$
which contradicts our uniform volume bound \eqref{volumebound}.
This contradiction yields the desired claim.

 From Theorem \ref{thmlocal1} and the Claim, there is a constant $C>0$,
independent of $T\in [0, T^*)$ such that $\forall (t,x)\in
[0,T]\times M$
\begin{equation}\label{eq6.7.2}
|\nabla u(t,x)|+ |\nabla^2 u(t,x)|\le C.
\end{equation}
Using the fact $\ds \int_M \s_1( g) dv( g)\equiv 1$, we have
$\forall (t,x)\in [0,T]\times M$
$$
|u(t,x)|+ |\nabla u(t,x)|+ |\nabla^2 u(t,x)|\le C.
$$
Therefore, we finish the proof of Theorem. \qed

\begin{rem}
Our perturbed equation is so chosen such that 
the argument in Corollary \ref{thmglobal1} works and the estimate in Lemma \ref{lastlem} hold.
\end{rem}

\begin{rem}\label{thmglobal2}
Under the same assumptions as in Remark
\ref{thmlocal2},
 there is a  constant $C$ depending only on  $g_0$  such that
\begin{equation}
\|u\|_{C^2(M)}\le C.
\end{equation}
\end{rem}

\section{Uniform parabolicity}

We prove in this Section that our flow (\ref{flow}) preserves the positivity of the
scalar curvature. 

\begin{pro}
\label{positivity} There is a constant $C_0>0$, independent of $T\in
[0, T^*)$  such that $\s_1(g(t))> C_0$ for any $t\in [0, T]$.
\end{pro}

\pr The proof is a modification of  the proof given in \cite{GuanWang2} and
\cite{GeWang1}, with more attention on $\nu_1$ and $\nu_2$, and their derivatives.
 Recall
\[
\ba{l} \vs \ds W=(w_{ij})=(\n ^2_{ij} u+u_iu_j-\frac {|\n u|^2}2
 (g_0)_{ij}+(S_{g_0})_{ij}), \\
\ds \nu=\nu_1(g)e^{-(4-\e)u}+\nu_2(g). \ea
\]
We define
$$
F:= \frac{\s_2(W)-\nu}{\s_1(W)}-\kappa e^{-2u}
$$
for some sufficiently large $\kappa$ to be fixed later. Hence, $
F=u_t -(\kappa+\mu(g))e^{-2u}-m(g(t))$. By Corollary
\ref{thmglobal1} one can show that there is a constant $c_1>0$ which
is independent of $T>0$ such that
\begin{equation}\label{eq_add2}
\frac 1{c_1}>\int_M(\s_2(g)-\e e^{4u})dv(g)>c_1, \quad \frac 1{c_1} >\int_M
e^{\e u} dv(g)-\e (\int_M \s_1(g) dv(g))^{3-\e}>c_1.
\end{equation}
To show this,  from  Corollary \ref{thmglobal1} we first have that $\int_M(\s_2(g)-\e
e^{4u})dv(g)$ and $\int_M e^{\e u} dv(g)-\e (\int_M \s_1(g)
dv(g))^{3-\e}$ are bounded from above by some positive constants. It follows from Proposition \ref{pro} that ${\cal
E}_\e(g)$ is bounded from below by some positive constant and the
fact that $\int \s_1(g)dv(g)$ is constant along the flow. Therefore,
the second part in the inequalities yields. As a consequence,
$\nu_1(g)$, $\nu_2(g)$ and $\mu(g)$ are bounded from above and from
below by some positive constants. Again from Corollary
\ref{thmglobal1}, $m(g)$ is bounded.

Without loss of generality, we assume that the minimum of
$F$ is achieved at $(x_0, t_0)\in M\times (0, T]$.
 Near $(x_0,t_0)$, we have
\beq \label{v1}\ba{rcl} \ds\vs \frac{d}{dt} F &=& \ds \sum_{ij}
A^{ij} (\n_g^2(u_t))_{ij}+ 2\kappa e^{-2u}u_t+\frac{(4-\e) \nu_1(g) e^{-(4-\e)u}u_t}{\s_1(W)}
-\frac {\alpha(t)}{\s_1(W)}\\
&= &\ds\vs \sum_{ij} A^{ij} \left[(\n_g^2(F))_{ij} +(\kappa+ \mu(g))
(\n_g^2(e^{-2u}))_{ij}\right]+ 2\kappa e^{-2u}u_t+\frac{(4-\e)
\nu_1(g) e^{-(4-\e)u}u_t}{\s_1(W)}\\
&& \ds -\frac {\alpha(t)}{\s_1(W)}, \ea \eeq where
$$
\alpha(t):=\frac{d\nu_1(g)}{dt}e^{-(4-\e)u}+\frac{d\nu_2(g)}{dt}
$$
and
\[
A^{ij}:=\frac{\p F}{\p
w_{ij}}=\frac{(\s_1^2(W)-\s_2(W)+\nu)\delta^{ij}- \s_1(W)W^{ij}}
{\s_1^2(W)}
\]
is positive definite.  We choose the normal coordinates so that
$W$ is a diagonal matrix at $(x_0,t_0)$. First we claim there exists
some constant $c_2>0$ independent of $T$ and $\kappa$ such that for
all $t\in[0,T]$
\begin{eqnarray}
\label{eq_add3} \left|\frac{d\nu_1(g)}{dt}(t)\right|\le
c_2(1+\kappa+\frac1{\s_1(W)(x_0,t_0)}),\\
\label{eq_add4} \left|\frac{d\nu_2(g)}{dt}(t)\right|\le
c_2(1+\kappa+\frac1{\s_1(W)(x_0,t_0)}).
\end{eqnarray}
Using \eqref{flow}, \eqref{eq_add2} and Corollary \ref{thmglobal1},
we can estimate
\beq\label{eq_add}
\left|\frac{d\nu_1(g)}{dt}(t)\right|\le c
(1+\int_M\frac1{\s_1(W)(x,t)} dv(g_0))\le c(1+\max_{x\in
M}\frac1{\s_1(W)(x,t)}).
\eeq
Since $(x_0,t_0)$ is the minimum of $F$ in $ M\times [0, T]$, we
have
\beq\label{add_x}
\frac{\nu}{\s_1(W)}(x,t)\le
\frac{\s_2(W)}{\s_1(W)}(x,t)-\kappa
e^{-2u(x,t)}-\frac{\s_2(W)-\nu}{\s_1(W)}(x_0,t_0)+\kappa
e^{-2u(x_0,t_0)},
\eeq
for any point $(x,t)\in M\times [0,T]$.
Applying Corollary \ref{thmglobal1}, we have that $\s_1(W)$, $\s_2(W)$
and  $e^{-2u}$ are bounded and $\nu$ is bounded from above and from
below by some positive constants. Together with the fact
$\s_2(W)(x,t)\le \frac 13 \s_1^2(W)(x,t)$, (\ref{add_x}) implies there exists $c_3>0$
independent of $T$ and $\kappa$ such that for all $(x,t)\in M\times
[0,T]$
\begin{equation}
\label{eq_add6} \frac1{\s_1(W)(x,t)}\le
c_3(1+\kappa+\frac1{\s_1(W)(x_0,t_0)}),
\end{equation}
which, in turn, together with (\ref{eq_add}), implies \eqref{eq_add3}. Similarly, we have \eqref{eq_add4}.
Hence, we prove the desired claim. As a consequence, we have at the
point $(x_0,t_0)$
\begin{equation}
\label{eq_add5} \frac {|\alpha(t_0)|}{\s_1(W)}\le c
(\frac{1+\kappa}{\s_1(W)}+\frac{1}{\s_1^2(W)})
\end{equation}
Since $(x_0,t_0)$ is the minimum of $F$ in $ M\times [0, T]$,
 at this point we  have
$ \ds\frac {d  F}{dt} \le 0$,
 $F_l =0$ $\forall l$
and $ ( F_{ij}) \text{ is non-negative definite.}$ Note that
\[ ( \n^2_g)_{ij}F = F_{ij}+ u_iF_j+u_j F_i-\sum_l u_l F_l \d_{ij}=F_{ij},
\]
at $(x_0,t_0)$, where  $F_{j}$ and $F_{ij}$ are the first and second
derivatives
 with respect to the back-ground metric $g_0$. From the positivity of $A$ and (\ref{v1}), we have
\begin{equation}
\label{5add1}\begin{array}{rcl}
0& \ge &  \vs\ds F_t- \sum_{i,j} A^{ij} F_{ij} \\
&&\\
 &\ge&\vs\ds (\kappa+ \mu(g))\sum_{i,j}  A^{ij}
\{(e^{-2u})_{ij}+u_i(e^{-2u})_{j}+u_j(e^{-2u})_{i}-\sum_l u_l(e^{-2u})_l\d_{ij}\}\\
&&\vs\ds + 2\kappa e^{-2u}u_t+\frac{(4-\e) \nu_1(g) e^{-(4-\e)u}u_t}{\s_1(W)}-\frac \alpha{\s_1(W)}\\
&=&\vs\ds (\kappa+ \mu(g)) e^{-2u}\sum_{i,j}  A^{ij}\{- 2w_{ij}+2u_iu_j+2S(g_0)_{ij}+|\n u|^2\d_{ij}\}\\
&&\vs\ds + 2\kappa e^{-2u}u_t+\frac{(4-\e) \nu_1(g) e^{-(4-\e)u}u_t}{\s_1(W)}-\frac \alpha{\s_1(W)}\\
&= & \ds \vs(\kappa+ \mu(g))
e^{-2u}\left(\frac{-2\s_2(W)-2\nu}{\s_1(W)}\right)  + 2\kappa e^{-2u}u_t+\frac{(4-\e) \nu_1(g) e^{-(4-\e)u}u_t}{\s_1(W)}\\
&&\vs\ds+(\kappa+ \mu(g)) e^{-2u} \sum_{i,j}
A^{ij}(2u_iu_j+2S(g_0)_{ij}+|\n u|^2\d_{ij})-\frac \alpha{\s_1(W)}.
\end{array}\end{equation}
Here we have used $\ds\sum_{i,j}  A^{ij}
w_{ij}=\frac{\s_2(W)+\nu}{\s_1(W)}$. A direct computation gives
\begin{equation}
\label{5add2}
\begin{array}{lll}
\ds\sum_{i,j} A^{ij}S(g_0)_{ij}&=&\ds\vs
\frac{(\s_1^2(W)-\s_2(W))\s_1(g_0)}{\s_1^2(W)}-
\frac1{\s_1(W)}\sum_{i,j}W^{ij}S(g_0)_{ij}
+\frac{\nu\s_1(g_0)}{\s_1^2(W)} .
\end{array}
\end{equation}

Gathering (\ref{eq_add5}), (\ref{5add1}) and (\ref{5add2}), we have
\begin{equation}
\label{5add3}\begin{array}{rcl}
0& \ge &  \vs\ds F_t- \sum_{i,j} A^{ij} F_{ij} \\
&\ge&\vs\ds (\kappa+ \mu(g))
e^{-2u}\left[\frac{-2\s_2(W)-2\nu}{\s_1(W)}
+\frac{2(\s_1^2(W)-\s_2(W))
\s_1(g_0)}{\s_1^2(W)}\right.\\
&&\left.\ds\vs- \frac2{\s_1(W)}\sum_{i,j}W^{ij}S(g_0)_{ij}
+\frac{2\nu\s_1(g_0)}{\s_1^2(W)} \right]\\
&&\ds -C(1+\frac{1+\kappa}{\s_1(W)}+\frac{1}{\s_1^2(W)}),
\end{array}\end{equation}
since $(A^{ij})$ is positive definite and $\kappa+ \mu(g)$ is
positive. Let us use $O(1)$ denote terms with a uniform  bound. One
can check again $\s_2(W)=O(1)$ for $\|u\|_{C^2}$ is uniformly
bounded and $\ds \sum_{i,j}W^{ij}S(g_0)_{ij}=O(1)$. Also the term
$\s_1^2(W)-\s_2(W)$ is always non-negative. We can choose $\kappa$
such that
$$
\frac{(\kappa+ \mu(g))\nu\s_1(g_0)e^{-2u}}{\s_1^2(W)} \ge
\frac{C}{\s_1^2(W)}
\quad\mbox{and
}\quad\kappa+ \mu(g)\ge \frac{\kappa +1}2
$$
Fixing such $\kappa$, from (\ref{5add3}) we conclude there holds at the point $(x_0,t_0)$
$$
0 \ge  \frac{\kappa}{\s_1^2(W)}-c_4(\frac{\kappa+1}{\s_1(W)}+1)
$$
for some positive constants $c_4>0$ independent of $T$ and $\kappa$. Consequently, there
is a positive constant $c_5>0$ (independent of $T$) such that
$$
\s_1(W)(x_0,t_0)\ge c_5.
$$
Hence, from \eqref{eq_add6} and Corollary \ref{thmglobal1},
there is a positive constant $c_6>0$, independent of $T$, such that
for all $(x,t)\in M\times [0,T]$
\[
\s_1(W)(x,t)\ge c_6.
\]
This finishes the proof of the Theorem. \qed


\section{Proof of main Theorems}

Now we can show the convergence of flow \eqref{flow}.
 \begin{thm}\label{existence} For small $\e>0$ flow \eqref{flow} with an initial metric
 $g_1\in {\cal C}_{1,\e}([g_0])$ converges to a
 metric $g_\infty$ satisfying \eqref{p-eq}.
 \end{thm}
 \pr With the $C^2$ estimates (Corollary \ref{thmglobal1}) and the uniform parabolicity
 (Proposition \ref{positivity}), one can show the convergence
 like in \cite{GuanWang2}.
 \qed

In order to estimate the value of $M_\e$ we need the following
\begin{lem}\label{lem4} There exists some $C_0>0$ depending only on $g_0$
such that for any $g=e^{-2u}g_0\in {\cal
C}_1([g_0])$ satisfying $\ds\int_M \s_2(g)dv(g)\ge0$ there holds
\begin{equation}
\label{equivnorm} C_0e^{\max u}\le \int e^{4u}dv(g)\le e^{\max
u}vol(g_0).
\end{equation}
\end{lem}

\noindent{\it Proof.} The second inequality is clear. We
prove the first inequality. As in \cite{GLW}, we have for all
$g\in {\cal C}_1([g_0])$ \beq
\begin{array}{rcl}
\label{eq6.9} \ds\int \s_2(g) dv(g) &\leq & \ds -\frac 1 {16}
\int |\n u|^4_{g_0} e^{4u} d(g)+c \int e^{4u}dv(g),
\end{array}
\eeq for some positive constant $c>0$. Since $ \ds\int \s_2(g)
dv(g)$ is non-negative,  we have 
\beq \label{eq6.10}
 4^4\int|\n e^{u/4}|_{g_0}^4dv(g_0)=\int |\n u|^4_{g_0} e^{4u} dv(g)\le c \int e^{4u}dv(g)= c \int (e^{u/4})^4dv(g_0),
\eeq which implies, with the help of  Sobolev's embedding Theorem ($W^{1,4}\subset C^{1/4}$),  for all
$x,y\in M$ \beq \label{eq6.11}
 |e^{u(x)/4}- e^{u(y)/4}| \le c (\int e^{4u}dv(g))^{1/4}(d_{g_0}(x,y))^{1/4},
\eeq where $d_{g_0}(x,y)$ is the distance between $x$ and $y$ with
respect to the metric $g_0$. Set \beq \label{eq6.13}
 \beta := e^{\max_M u}=  e^{u(x_0)}
\eeq for some $ x_0\in M$. It follows from (\ref{eq6.11}) that there
exists some $r>0$ independent of $u$ such that for any $y\in B( x_0,r)$
 \beq
\label{eq6.14} e^{u(y)/4}\ge \frac12 \beta^{1/4} \eeq 
Here the
geodesic ball $B( x_0,r)$ is taken for the metric $g_0$.
Hence, we
deduce
$$
\int_M e^{4u}dv(g)\ge \int_{B(x_0,r)} e^{4u}dv(g)\ge c\beta.
$$
Therefore, we have finished to prove the Lemma.\qed

Now we estimate the value of $M_\e$. The proof likes one given for Lemma \ref{lem2}, with the help of
 Lemma \ref{lem4}.
 \begin{lem}\label{lastlem}  Let  $C_0>0$ be the constant given in Lemma  \ref{lem4}. Any solution $g\in {\cal C}_{1,\e}$ of (\ref{p-eq}) satisfies \eqref{est-add}, i.e.,
 $$
{\cal E}_\e(g) \le (\frac 2{C_0\e})^\e\frac1{3(1-\e)}.$$
 \end{lem}

\pr 
 Multiplying
(\ref{p-eq}) by $e^{\e u}$ and integrating over
$M$, we have
\begin{equation}
\int_M \frac {e^{\e u}\s_2(g)}{\s_1(g)}dv(g)+\mu(g)\int_M e^{\e
u}dv(g)=\int_M \frac{(\nu_1(g)e^{2\e u}+\nu_2(g)e^{(4+\e)
u})}{\s_1(g)}dv(g)
\end{equation}
By the Cauchy-Schwarz inequality, we have
$$
\begin{array}{lll}
&\ds\int_M \frac{(\nu_1(g)e^{2\e u}+\nu_2(g)e^{(4+\e)
u})}{\s_1(g)}dv(g)\int_M {\s_1(g)}dv(g)
\\
\ge  &\ds\int_M \frac{\nu_1(g)e^{2\e u}}{\s_1(g)}dv(g) \int_M
{\s_1(g)}dv(g)\\
\ge &\ds\nu_1(g)(\int_M e^{\e u}dv(g))^2.
\end{array}$$
The above two inequalities implies that
\begin{equation}\label{z1}
\begin{array}{rcl}
\ds \int_M \frac {e^{\e u}\s_2(g)}{\s_1(g)}dv(g) &\ge&\ds\vs\nu_1(g)\frac{(\int_M e^{\e u}dv(g))^2}
{\int_M\s_1(g)dv(g)} -\mu(g)\int_M e^{\e
u}dv(g)\\
&=&\ds (1-\e)\frac {k(g)\int_M (\s_2(g)-\e e^{4u}) dv(g) \int_M e^{\e u} dv(g)}{\int_M \s_1(g)dv(g)}.
\end{array}
\end{equation}
In the last equality we have used the definitions of $\nu_1(g)$ and $\mu(g)$.

On the other hand, we recall the facts
$$
\s_2(g)\le \frac13(\s_1(g))^2
$$
and for all $g\in {\cal C}_{1,\e}([g_0])$
$$
\e \int_M e^{4u }dv(g)\le \int_M \s_2(g) dv(g).
$$
Hence, we get from Lemma \ref{lem4}
\beq\label{z2}
\begin{array}{lll}
\ds\int_M \frac{e^{\e u}\s_2(g)}{\s_1(g)}dv(g)&\ds\le \frac13 \int_M
e^{\e u}\s_1(g)dv(g)\le \frac13 e^{\e \max u}\int_M \s_1(g)dv(g)\\
\ds &\ds\le \frac13 (\int_M (\s_2(g)-\frac\e2 e^{4u}dv(g))^\e (\frac 2{C_0\e})^\e
\int_M \s_1(g)dv(g).
\end{array}
\eeq
\eqref{z1} and \eqref{z2} give us
\[
\frac {\int_M (\s_2(g)-\e e^{4u})dv(g) \int_M e^{\e u} dv(g)}
{(\int_M (\s_2(g)-\frac \e 2 e^{4u}) dv(g))^\e (\int_M \s_1(g) dv(g))^2} \le (\frac 2{C_0\e})^\e\frac1{3(1-\e)k(g)},
\]
which implies
$$
{\cal E}_\e(g)\le (\frac 2{C_0\e})^\e\frac1{3(1-\e)k(g)} \le  (\frac 2{C_0\e})^\e\frac1{3(1-\e)},
$$
since $k(g)\ge 1$.
 This yields the desired result.  \qed



\

\noindent{\it Proof of Theorem \ref{mainthm2}.}
If  any $g\in {\cal C}_{1}([g_0])$ satisfies $\ds\int_M \s_2(g)\le 0$,
then $\widetilde Y([g_0])\le 0<1/3$. Hence we consider that there is $g\in {\cal
C}_{1}([g_0])$ with $\ds\int_M \s_2(g) > 0$.
For such a metric $g$ we can choose a small number $\e_0>0$ such that 
 $g\in {\cal C}_{1,\e}([g_0])$ for any $\e \in (0,\e_0)$. Hence, we have
 \[
 {\cal E}_\e(g) \le M_\e.\]
 Theorem \ref{existence} and Remark 3 imply that $M_\e$ is achieved by a metric $\tilde g \in {\cal C_1}\cap [g_0]$ satisfying (\ref{p-eq}).
 From Lemma  \ref{lastlem} we have \[
M_\e={\cal E}_\e(\tilde g)\le (\frac 2{C_0\e})^\e\frac1{3(1-\e)},
\]
and hence  
 \[{\cal E}(g)\le \frac 13.\]
 Therefore we have 
 \[\widetilde Y([g_0]) \le \frac 1 3.\]
 This finishes the proof of the Theorem. \qed

We consider now the energy functional ${\cal E}$  in a larger
class $$\overline{{\cal C}_1([g_0])}:=\{g=e^{-2u}g_0\,|\, R\ge
0\}$$
and define
\begin{equation}
\ds \bar Y([g_0]):=\sup_{ g\in \overline{{\cal C}_1([g_0])}}\frac
{\ds vol(g)\int_M \s_2( g) dv( g)} {\ds (\int_M \s_1( g) dv( g))^2}.
\end{equation}
Note that in $\overline{{\cal C}_1([g_0])}$ there is no metric with $R\equiv 0$, if $g_0\in {\cal C}_1$.
We have the following result, which improves slightly Theorem
\ref{mainthm2}.

\begin{thm} \label{mainthm2bis} If $(M^3, g_0)$ is a closed $3$-dimensional manifold with positive
Yamabe constant $Y_1([g_0])> 0$, then
\begin{equation}\label{Invariant2}
\ds \bar Y([g_0])\le \frac 1 3.
\end{equation}
Moreover, equality holds if and only if $(M^3, g_0)$ is space form.
\end{thm}

 \

\noindent {\it Proof of Theorem \ref{mainthm2bis}.} For any metric
$g=e^{-2u}g_0\in \overline{{\cal C}_1([g_0])}$, we consider
$g_t=e^{-2tu}g_0$ for $0<t<1$. Clearly, $g_t\in {\cal C}_1([g_0])$.
By the approximation arguments and Theorem \ref{mainthm2}, we have
$$
{\cal E}(g)=\lim_{t\to 1}{\cal E}(g_t)\le \frac13.
$$
Now we suppose ${\cal E}(g)= \frac 13$. Thus, $g$ is an extremal
metric in the class of $\overline{{\cal C}_1([g])}$ for the energy
functional ${\cal E}$. Denote $M_1:=\{x\in M, \s_1(g)(x)=0\}$ and
$M_2:=\{x\in M, \s_1(g)(x)>0\}$. We have $M=M_1\cup M_2$ and
(\ref{new}) is verified in $M_2$. On the other hand, if $x\in M_1$,
then $\s_1(g)(x)=0$ and  $\s_2(g)(x)\le 0$. Hence, we deduce in
$M_1$
$$
\s_2(g)- 3r_2(g)+2s(g)\s_1(g)< 0
$$
since $r_2(g)> 0$. On the other hand, by the definition of $r_2$ and $s$ we know
$$
\int_M (\s_2(g)- 3r_2(g)+2s(g)\s_1(g))dv(g)= 0
$$
which implies  $M=M_2$ and we have Equation (\ref{new}). Therefore,
from Lemma \ref{lem2}, we infer that $M$ is an Einstein manifold and we
finish the proof. \qed

\

\noindent {\it Proof of Theorem \ref{mainthm1}.} 
Let $(M^3,g)$ be a metric of non-negative scalar curvature.
 If its Yamabe constant is positive, then by
Theorem \ref{mainthm2bis} we have ${\cal E}(g)\le \bar Y([g]) \le
1/3$, which is equivalent to \eqref{eq1} by Lemma \ref{keylemma}. Hence
we only need to consider the case that $g$ has zero Yamabe
constant. In this case one can show that $g$ has scalar curvature
zero and hence \eqref{eq1} holds trivially.

It is trivial to see that an Einstein metric satisfies \eqref{eq1}
with equality. Now assume that $g$ is a metric of non-negative
scalar curvature which satisfies \eqref{eq1} with equality. If
$\int_M \s_1(g) dv(g)=0$, then we have $\s_1(g)=0$, which implies
that $\bar R=R\equiv 0$.  Using \eqref{eq1},  $g$ is a Ricci flat
metric, and hence a flat metric. If $\int_M \s_1(g) dv(g)>0$, by Lemma \ref{keylemma} we have
${\cal E}(g)=1/3$. Hence $(M^3,g) $ is a space form   by Theorem
\ref{mainthm2bis}.
\qed

In a recent joint work with Xia \cite{GWX} we proved  the rigidity of \eqref{eq1}, namely under the conditions in Theorem A
equality in \eqref{eq1} holds if and only if $(M,g)$ is an Einstein metric.

\

\section{Problems and
Conjectures}

We end the paper by proposing several related  problems and  conjectures. 

\

\noindent{\it Conjecture 1.} {\it Theorem \ref{mainthm1} holds if $g$ has a non-negative first Yamabe constant $Y_1([g])$}.

\

When $g$ has a negative first Yamabe constant, Theorem \ref{mainthm1} is not true. For example see \cite{DT}.

\noindent{\it Problem 1.} {\it  $\widetilde Y([g_0])$  is achieved.}

\

This is a Yamabe type problem, but with a different property.
From the analysis developed here, together with a classification
result of blow-up solutions like in \cite{LL2}, one can expect  that
this conjecture is true if $\widetilde Y([g_0])<1/3$.
It is trivial to see that any metric $g$ with constant sectional curvature  satisfies ${\cal E}(g)=1/3$, and hence
\[ \widetilde Y([g])=\frac 13= \widetilde Y([g_{\S^3}]),\]
where $g_{\S^3}$ is the standard round metric on $\S^3$. We conjecture 

\

\noindent{\it Conjecture 2.} {\it Let $(M^3,g_0)$ be a closed manifold with $g_0\in {\cal C}_1$. If $\tilde Y([g_0])=1/3$,
then $(M^3,g_0)$ is conformally equivalent to a  3-dimensional spherical space form.
}

\

It is inetresting to see that this conjecture, if it is true, gives a characterization of a conformal Einstein metric on a 3-dimensional manifold.

\

Let
\[J(g):=\int_M\s_1(g) dv(g) \cdot \int_M
\s_2(g) dv(g).\]

\

\noindent{\it Conjecture 3.} {\it Let $(M^3,g_0)$ be a closed manifold with $g_0\in {\cal C}_1$. The  following statement 
\[\widetilde Y_{2,1}([g_0]):=\sup_{g\in {\cal C}_1([g_0])} J(g) \le J(g_{\S^3})\]
is true.
}

\


This conjecture is closely related to a problem which was asked by Viaclovsky to us several years ago.
Let
\[J_2(g):=vol(g)^{1/3}\cdot{ \int_M
\s_2(g) dv(g)}
.\]
He asked if $\sup_{g\in {\cal C}_2([g_0])} J_2(g)$ is bounded. We believe that it is true
and we even believe more.

\

\noindent{\it Conjecture 4.} {\it Let $(M^3,g_0)$ be a closed manifold with $g_0\in {\cal C}_1$. The following statement 
\[\widetilde Y_{2}([g_0]):=\sup_{g\in {\cal C}_1([g_0])} J_2(g) \le J_2(g_{\S^3})\]
is true.
}

\

It is easy to see that Conjecture 3 implies  Conjecture 4 and Conjecture 4 implies 
Theorem \ref{mainthm2}. 

The Euler-Lagrange equation of $J_2$ is the so-called $\s_2$-Yamabe equation
\begin{equation}
 \label{sigma2}
\s_2(g)=b,
\end{equation}
for soma constant $b$. A Lemma  \ref{lem2} type result is true for this equation. This  in fact directly follows from
a volume comparison result of Gursky-Viaclovsky \cite{GuV3}, which in turn follows from a volume comparison Theorem of Bray \cite{Bray}.

\begin{lem} Let $g\in {\cal C}_1$ be a metric on a 3-dimensional manifold $M^3$ satisfying (\ref{sigma2}) then
 \[  J_2(g) \le J_2(g_{\S^3}).\]
\end{lem}
\pr 
We need only to consider the case $b>0$, otherwise the Lemma is trivial.
We may assume that $b= \s_2(g_{\S^3})$, i.e, $\s_2(g)=\s_2(g_{\S^3})$. Theorem 1.2 in \cite{GuV3} implies
\[vol(g)\le vol (g_{\S^3}).\]
Hence we have
\[ J_2(g) \le (vol(g_{\S^3}))^{\frac 43} \ \s_2(g_{\S^3})=J_2(g_{\S^3}).
 \]
\qed
Hence, to show Conjecture 4, as inspired by the proof given above, one needs only either to show that $J_2$ is achieved, or
to show a suitable perturbed functional has a maximum, together with a Lemma \ref{lastlem} type estimate. This is a difficult problem.
There is even an extra difficulty that the corresponding flow is in general not parabolic. However for functional $J$ there is no this extra difficulty.
This is known from  Lemma \ref{lem3.1}. Therefore, it may be better to study Conjecture 3 first. 

The Euler-Lagrange equation of $J$ is the so-called quotient  equation
\begin{equation}
 \label{quo}
\frac{\s_2(g)}{\s_1(g)}=b.
\end{equation}
With the same idea, we  need  the following comparison result.

\

\noindent{\it Conjecture 5.} {\it Let $(M^3,g)$ be a closed 3-dimensional manifold with $g\in {\cal C}_1$. 
Assume that
\[
 \frac{\s_2(g)}{\s_1(g)}\ge  \frac{\s_2(g_{\S^3})}{\s_1(g_{\S^3})}.
\]
Then
\[
 \int_M \s_1(g)  dv(g) \le \int_{\S^3} \s_1(g_{\S^3}) dv(g_{\S^3}).
\]

}

\


\

If these conjectures are true, then
it is  natural to ask

\

\noindent{\it Problem 2.} {\it Are   $\widetilde Y_{2,1}([g_0])$  and $\widetilde Y_{2}([g_0])$ achieved?}

\


\begin{thebibliography}{99}

\bibitem{Bray}
H.L. Bray, The Penrose inequality in general relativity and volume comparison theorems involving scalar curvature, Dissertation, Stanford University, 1997.

\bibitem{CGY} A. Chang, M. Gursky and P. Yang, An equation of Monge-amp\`ere type in conformal geometry, and four manifolds of positive Ricci curvature, Ann. of Math. \textbf{155} (2002) 709--787
\bibitem{CGY2} A. Chang, M. Gursky and P. Yang, A conformally invariant sphere theorem in four dimensions, Publ. Math. Inst. Hautes \'Etudes Sci. \textbf{98} (2003) 105--143
\bibitem{Chen} S. Chen, Local estimates for some fully nonlinear elliptic equations, Int. Math. Res. Not.
 \textbf{2005} (2005), 3403--3425,

\bibitem{CD}  G. Catino,  and Z. Djadli,
 Conformal deformations of integral pinched 3-manifolds.  Adv. Math. \textbf{223}  (2010), 393--404
\bibitem{DT}C. De Lellis and P. Topping,
Almost Schur Theorem, \textbf{Arxiv 1003.3527}.

\bibitem{GLW} Y. Ge, C.-S. Lin and G. Wang, On $\s_2$-scalar curvature, J. Diff. Geom., \textbf{84}  (2010),  45--86.

\bibitem{GeWang2} Y. Ge and G. Wang, On $\s_2$-scalar curvature II, preprint (2007).
\bibitem{GeWang1} Y. Ge and G. Wang,
 On a conformal quotient equation.  Int. Math. Res. Not. IMRN  2007,   Art, ID rnm019, 32 pp.
\bibitem{GeWang_Proc} Y. Ge and G. Wang, An almost Schur Theorem on 4-dimensional manifolds, preprint (2010).
\bibitem{GWX} Y. Ge, G. Wang and Chao Xia, On  problems related to an inequality of De  Lellis and Topping,
preprint.


\bibitem{GuanLinWang} P. Guan, C.-S. Lin and G. Wang,
 Application of The Method of Moving Planes to Conformally Invariant Equations, Math. Z. \textbf{247} (2004) 1--19

\bibitem{guan}P. Guan, Topics in Geometric Fully Nonlinear Equations, Lecture Notes, http://www.math.mcgill.ca/guan/notes.html
\bibitem{GuanWang1} P. Guan and G. Wang,
 Local estimates for a class of fully nonlinear equations arising from conformal geometry,
Int. Math. Res. Not. \textbf{2003} (2003), 1413--1432,
\bibitem{GuanWang2} P. Guan and G. Wang, A fully nonlinear conformal flow on locally conformally flat manifolds, J. Reine Angew. Math., \textbf{557} (2003) 219--238


\bibitem{Gu} M. Gursky, The principal eigenvalue of a conformally invariant differential operator, with an application to semilinear elliptic PDE, Comm. Math. Phys. {\bf 207} (1999), 131--143.

\bibitem{Gu2} M. Gursky, Fully nonlinear equations,
ellipticity, and curvature pinching,
S\'eminaires $\&$ Congr\`es, {\bf 19} (2008) 31--45


\bibitem{GuV2} M. Gursky and J. Viaclovsky,  A new variational characterization of three-dimensional space forms,
  Invent. Math., {\bf 145}  (2001),   251--278. 

\bibitem{GuV} M. Gursky and J. Viaclovsky, 
 A fully nonlinear equation on four-manifolds with positive scalar curvature.  J. Diff. Geom.  {\bf 63}  (2003), 131–154.

\bibitem{GuV3} M. Gursky and 
J. Viaclovsky, Volume comparison and the $\sigma_k$-Yamabe problem, Advances in Math. {\bf 187} (2004) 447--487.

\bibitem{LL} A. Li and Y. Li, On some conformally invariant fully nonlinear equations, Comm. Pure Appl. Math.,
\textbf{56} (2003) 1416--1464

\bibitem{LL2} A. Li and Y. Li, On some conformally invariant fully nonlinear equations. II. Liouville, Harnack and Yamabe, Acta Math., \textbf{195} (2005) 117--154
\bibitem{V}J. Viaclovsky, Conformal geometry, contact geometry, and the calculus of variations, Duke Math. J., {\bf 101} (2000), 283--316.

\bibitem{Survey} J. Viaclovsky, Conformal geometry and fully nonlinear equations, Inspired by S. S. Chern, 435--460, Nankai Tracts Math. 11 World Sci. Publ., Hackensack, NJ, 2006
\bibitem{Wang}
G. Wang, $\sigma_k$-scalar curvature and eigenvalues of the Dirac operator, Ann. Global Anal. Geom. \textbf{30}
 (2006) 65--71

\end{thebibliography}
\end{document}